\documentclass[12pt]{amsart}
\usepackage{amscd,amssymb}
\usepackage[arrow,matrix]{xy}
\usepackage[plainpages,backref,urlcolor=blue]{hyperref}

\theoremstyle{plain}
\newtheorem{thm}[subsection]{Theorem}

\theoremstyle{definition}
\newtheorem{rk}[subsection]{Remark}
\newtheorem{definition}[subsection]{Definition}
\newtheorem{ex}[subsection]{Example}

\numberwithin{equation}{section}
\newcommand{\A}{{\mathcal A}}

\newcommand{\C}{\mathbb{C}}

\newcommand{\PP}{\mathbb{P}}

\newcommand{\wI}{\widehat{I}}
\newcommand{\wJ}{\widehat{J}}

\begin{document} 

\title [Free divisors versus stability and coincidence thresholds]
{ Free divisors versus stability and coincidence thresholds }

\author[Gabriel Sticlaru]{Gabriel Sticlaru}
\address{Faculty of Mathematics and Informatics,
Ovidius University, Constanta, Romania}
\email{gabrielsticlaru@yahoo.com }



\begin{abstract} We show that there is an unexpected relation between free divisors and 
stability and coincidence thresholds for projective hypersurfaces.

\end{abstract}

\maketitle
 
\section{Introduction}

Let $S=\C[x_0,...,x_n]$ be the graded ring of polynomials in $x_0,,...,x_n$ with complex coefficients and denote by $S_r$ the vector space of homogeneous polynomials in $S$ of degree $r$. 
For any polynomial $f \in S_r$ we define the {\it Jacobian ideal} $J_f \subset S$ as the ideal spanned by the partial derivatives $f_0,...,f_n$ of $f$ with respect to $x_0,...,x_n$. For $n=2$ we use $x,y,z$ instead of
$x_0, x_1, x_2$ and $f_x, f_y, f_z$  instead of $f_0, f_1, f_2$.

The Hilbert-Poincar\'{e} series of a graded $S$-module $M$ of finite type is defined by 
\begin{equation} 
\label{eq2}
HP(M)(t)= \sum_{k\geq 0} \dim M_k\cdot t^k 
\end{equation} 
and it is known, to be a rational function of the form 
\begin{equation} 
\label{eq3}
HP(M)(t)=\frac{P(M)(t)}{(1-t)^{n+1}}=\frac{Q(M)(t)}{(1-t)^{d}}.
\end{equation}

For any polynomial $f \in S_r$ we define the corresponding graded {\it Milnor} (or {\it Jacobian}) {\it algebra} by
\begin{equation} 
\label{eq1}
M=M(f)=S/J_f.
\end{equation}

Smooth hypersurfaces of degree $d$ have all the same  Hilbert-Poincar\'{e} series, namely
\begin{equation}
HP(M(f))(t)=\frac{(1-t^{d-1})^{n+1}}{(1-t)^{n+1}}=(1+t+t^2+\ldots + t^{d-2} )^{n+1}.
\end{equation} 
As soon as the hypersurface $V(f)$ acquires some singularities, the series $HP(M(f))$ is an infinite sum.

Let $ \tau(V(f))=\sum_{j=1,p}\tau(V(f),a_j) $
be the global Tjurina number of the hypersurface $V(f)$. In particular, the
Hilbert polynomial $H(M(f))$ is constant and this constant is $\tau(V(f))$.

\bigskip
For a hypersurface $D: f=0$ in $\PP^n$ with isolated singularities we recall the following invariants, introduced in 
 \cite{DS2}.

\begin{definition}
\label{def}

\noindent (i) The {\it coincidence threshold} $ct(D)$ defined as
$$ct(D)=\max \{q~~:~~\dim M(f)_k=\dim M(f_s)_k \text{ for all } k \leq q\},$$
with $f_s$  a homogeneous polynomial in $S$ of degree $d=\deg f$ such that $D_s:f_s=0$ is a smooth hypersurface in 
$ \PP^n$.

\bigskip

\noindent (ii) The {\it stability threshold} $st(D)$ defined as
$$st(D)=\min \{q~~:~~\dim M(f)_k=\tau(D) \text{ for all } k \geq q\}$$
where $\tau(D)$ is the total Tjurina number of $D$, i.e. the sum of all the Tjurina numbers of the singularities of $D$.

\end{definition} 

\noindent 

We refer for the definiton of free divisors and their main properties to \cite{ST} and the extensive reference list there.

\section{Main result} 

The main result of this note is the following.

\begin{thm} \label{thm1}
Let $f \in S_d$ and set as usual $T=(n+1)(d-2)$. Assume that $st(D)<ct(D)$. Then the following hold.

\begin{enumerate}
\item The degree $d$ is odd, $ ct(D)=(T+1)/2$ and $st(D)=(T-1)/2$.

\item If in addition $n=2$, then the curve $D:f=0$ is a free divisor.
\end{enumerate}

\end{thm}

To prove the first claim is fairly easy. Recall that the coefficients $c_k$ of the Hilbert-Poincar\'e series of
$M(f)$ form a (strict) log-concave sequence by  \cite{St1}, symmetric with respect to the middle point $T/2$, i.e. $c_k=c_{T-k}$. It follows that one has strict inequalities $c_{k-1}<c_k$ for all $k \leq T/2$. This implies that $st(D) \geq T/2$ when $d$ is even and $st(D) \geq (T-1)/2$ when $d$ is even. 

Suppose we are in the case $d$ is even. Then $st(D) > T/2$ implies $st(D) \geq ct(D)$ as beyond
$T/2$ the sequence $c_k$ is strictly decreasing. Hence one should have $st(D) = T/2$, which implies $ct(D)=st(D)$.

Suppose now that we are in the case $d$ is odd. Then again $st(D) > (T-1)/2$ implies $st(D) \geq ct(D)$ as beyond
$T/2$ the sequence $c_k$ is strictly decreasing. Hence one should have $st(D) = (T-1)/2$, which implies $ct(D)=st(D)+1= (T+1)/2 $.

The proof of the second claim is much more subtle.

Recall that $D$ is a free divisor if and only if the Jacobian ideal $J_f$ is a perfect ideal, see for instance the bottom of page 1 in \cite{ST}. Next, the Corollary 1.2 in  \cite{ST} implies the following fact: if the Castelnuovo-Mumford regularity $reg(M(f))$ satisfies
\begin{equation} 
\label{CM}
reg(M(f)) \leq (3d-7)/2,
\end{equation}
then $J_f$ is perfect and one has $reg(M(f)) \leq (3d-7)/2$ with $d$ odd.

Using \cite{D3}, Corollary 4, we know that
$$reg(M(f))=\max(T-ct(D),sat(J_f)-1),$$
where 
for any homogeneous ideal $I$, we let  $\wI$ be the associated saturated ideal and  we consider the graded artinian $S$-module
\begin{equation} 
\label{eqSD}
SD(I)=  \frac{\wI}{I},
\end{equation}
called the {\it saturation defect module} of $I$ and the {\it saturation threshold} $sat(I)$ defined as
\begin{equation} 
\label{eqsat}
sat(I)=\min \{q~~:~~\dim I_k= \dim \wI _k \text{ for all } k \geq q\}.
\end{equation}
Now, it follows from (1) above that $T-ct(D)=(T-1)/2=(3d-7)/2$, so we need to control $sat(J_f)$.
Corollary 2 in \cite{D3}  tells us that
$$sat(J_f)\leq max(T-ct(D),st(D)).$$
and hence in our situation $sat(J_f)-1\leq (T+1)/2-1=(T-1)/2= (3d-7)/2 $. Hence $D$ is a free divisor.
This completes the proof of Theorem \ref{thm1}.

\begin{rk}
\label{sat}
Note that in fact $D$ is a free divisor if and only if $J_f=\wJ_f$, see top of page 6 in \cite{ST}, and hence if and only if $sat(J_f)=0$. In particular, for a free divisor one has
$$reg(M(f))=T-ct(D).$$

\end{rk}

\section{Examples of hypersurfaces} 

We show here examples of hypersurfaces with $st(D) \leq ct(D)$

\begin{ex} (a triangle and a conic+tangent)
\label{ex1}
The example $D:f=xyz=0$ has $sat(J_f)=0$, $2=ct(D)>st(D)=1$. The conic plus tangent $D:f=x(xz+y^2)$ has the same invariants and Hilbert-Poincar\'{e} series $HP(M(f))(t)=1+3(t+\cdots .$
\end{ex}

\begin{ex} (two free line arrangements in $\PP^2$)
\label{ex1}
Consider the line arrangement
$$D(\A_3): f=(x^2-y^2)(y^2-z^2)(x^2-z^2)=0.$$
The Hilbert-Poincar\'{e} series is: 
$$HP(M(f))(t)=1+3t+6t^{2}+10t^{3}+15t^{4}+18t^{5}+19(t^{6}+\cdots .$$
Then $ct( D(\A_3))=st(D(\A_3))=6$ and the curve $D(\A_3)$ is free, since one can verify by Singular the vanishing $sat(J_f)=0$.
Similarly, for the line arrangement 
$$D: g= xyzf=xyz(x^2-y^2)(y^2-z^2)(x^2-z^2)=0,$$
one has $ct(D)=10 < st(D)=11$ and the Hilbert-Poincar\'{e} series is: $1+3t+6t^{2}+10t^{3}+15t^{4}+21t^{5}+28t^{6}+36t^{7}+42t^{8}+46t^{9}+48t^{10}+49(t^{11}+\cdots .$ So this curve is free for the same reason as above.
\end{ex}

\begin{ex} (a sequence of free divisors constructed in  \cite{ST}, Prop. 2.2.)
\label{ex2}
 Consider the sequence of free divisors $D_d: f=y^{d-1}z +x^d+x^2y^{d-2}=0$, for $ d\geq 5$
Then, for $d=5$ we have $st(D)=4<5=ct(D)$, hence we are in the setting of Theorem \ref{thm1} (ii).
For $d=6$, one has $st(D)=6=ct(D)$, while for $d>6$ one has $st(D)>ct(D)$. So Theorem \ref{thm1} (2) has no converse.
\end{ex}

\begin{ex} (Simis' free irreducible sextic, see \cite{Si})
\label{ex3}
Consider the sextic $D:f=0$, where
$$f=4(x^2+y^2+xz)^3-27(x^2+y^2)^2z^2.$$
Then $D$ has three $E_6$ singularities and one node $A_1$, hence
 $\tau(D)=19$.  Moreover one has $ct(D)=st(D)=6$ since the Hilbert-Poincar\'{e} series of $M(f)$ can be computed using for instance Singular and we get
$$HP(M(f))(t)=1+3t+6t^{2}+10t^{3}+15t^{4}+18t^{5}+19(t^{6}+\cdots .$$
This is one of the very few irreducible free divisors in $\PP^2$ known to this day.
\end{ex}

\begin{ex} (a sequence of curves $C_d$ with $st(C_d)= ct(C_d)$.)
\label{ex4}
Let $C_d:f_d=0$ with $f_d(x,y,z)=f=x^2y^2z^{d-4}+x^5z^{d-5}+y^5z^{d-5}+x^d+y^d$ for $d\geq 5$.
Then, for any  $d\geq 5$, the curve  $C_d$ has a unique singularity at the point $(0:0:1)$ of type $T_{2,5,5}$ with $\tau=10.$ Then a computation by Singular suggests that $ct(D_d)=st(D_d)=3d-9$ (verified for $5 \leq d \leq 15$). Moreover all these divisors are not free since $sat(J_f)$ is not zero (verified for $5 \leq d \leq 15$).

\end{ex}

\begin{ex} (a nodal hypersurface with one node)
\label{ex5}
If $D$ is a nodal hypersurface in $\PP^n$ having exactly one node, it follows from Example 4.3 (i) in \cite{DS3} that $st(D)= ct(D)=T$.
For nodal curve of degree $5$, $C: f=x^4y+x^3y^2+y^5+xy^2z^2+(x^2+xy+y^2)z^3=0$ our invariants are $ct(C)=st(C)=9$ and
the Hilbert-Poincar\'{e} series of $M(f)$ is
$$HP(M(f))(t)=1+3t+6t^2+10t^3+12t^4+12t^5+10t^6+6t^7+3t^8+(t^9+\cdots .$$

\end{ex}


\section{Conclusion} 

For a hypersurface $D: f=0$ in $\PP^n$ with isolated singularities, we compare the Hilbert-Poincar\'{e} series of 
the corresponding graded {\it Milnor} {\it algebra} with the smooth case. 
The series coincides up to a point $ct(D)$ and Choudary-Dimca Theorem \cite {CD} ensures that from an index $st(D)$ we have stability.

The main result is Theorem \ref{thm1}, which establishes a relation between free divisors and 
stability $st(D)$ and coincidence thresholds $ct(D)$.

In the vast majority of cases one has $ct(D) \leq st(D)$, see for instance the examples given above for the case of equality.
The opposite inequality $st(D)<ct(D)$ may occur but it is extremely rare. 

Recall that one of the main question in this theory is the construction of irreducible free divisor of arbitrarily high degree. Simis' celebrated example in degree $d=6$ is discussed in Example \ref{ex3}. The main conclusion of the given examples is that Theorem \ref{thm1} (2) has no converse, even in the class of free line arrangements. Moreover, the examples show that there are a lot of free and not free divisors satisfying the equality $st(D)=ct(D)$.

All the computations for this paper were made by Singular package, see  \cite{DGP} and \cite{GP}.

\bigskip
\noindent \Large\textbf{Acknowledgment}\\[2mm] 
\footnotesize The author is grateful to A. Dimca for suggesting this problem. \\[3mm]

\end{document}